\magnification=\magstep1     
\tolerance 10000 
\baselineskip=12pt
\def\Hyper#1{\hyper {\eskip #1}}
\def\eskip{\hskip.25em\relax}
\def\real{{\tt I\kern-.2em{R}}}
\def\nat{{\tt I\kern-.2em{N}}}
\def\hyper#1{\ ^*\kern-.2em{#1}}

\def\hypernat{{^*{\nat }}}

\def\Hyper#1{\hyper {\eskip #1}}
\def\power#1{{{\cal P}(#1)}}
\def\qed{{\vrule height6pt width3pt depth2pt}\par\medskip}
\def\pars{\par\smallskip}

\def\r#1{{\rm #1}}
\def\b#1{{\bf #1}}

\def\sig{{^\sigma}}
\hoffset=.25in
\hsize 6.00 true in
\vsize 8.85 true in
\font\eightrm=cmr9
\def\m@th{\mathsurround=0pt}

\def\id{\par\hangindent2\parindent\textindent}
\def\textindent#1{\indent\llap{#1}}
\centerline{\bf Nonstandard Consequence Operators}\par
\centerline{\bf Generated By Mixed Logic-Systems}\par\bigskip 
\centerline{Robert A. Herrmann}\par\medskip
\centerline{Mathematics Department}
\centerline{U. S. Naval Academy}
\centerline{572C Holloway Rd.}
\centerline{Annapolis, MD 21402-5002}
\centerline{14 JUNE 2004. Latest revision 20 AUG 2010}\bigskip
{\leftskip=0.5in \rightskip=0.5in \noindent {\eightrm {\it Abstract:} In this paper, mixed logic-systems are modeled by finite collections of n-tuples and contain both standard and nonstandard coordinate values. It is shown that each of the specifically defined mixed logic-systems generates an internal nonstandard consequence operator. Hyperfinite mixed logic-systems are also investigated.\par}}\par\bigskip
 \noindent {\bf 1. Introduction.}\par\medskip  Abraham Robinson is the first to apply nonstandard analysis to the language of the ``Lower predicate calculus'' (Robinson, 1963, p. 84, Geiser, 1968). In this paper as well as others mentioned in the references, nonstandard analysis is applied to universal logics and informal languages. In particular, nonstandard methods are applied to finite consequence operators. All of the definitions not presented in this paper can be found in Herrmann (2006, 2001, 1993, 1991). In Herrmann (2006, 2001), logic-systems are defined for a nonempty specific informal language $\r L.$ These logic-systems are only composed of members of $\r L.$ A need has now arisen that requires the logic-system notion to be extended so that, after the language is properly embedded into the nonstandard structure - the Grundlegend or Extended Grundlegend Structure (Herrmann, 1993), then intuitively the logic-system contains elements from $\Hyper {\b L} - {\b L}.$ (Recall, that by the identification process $\sig {\b L} = \b L.$) \par\medskip 

\noindent {\bf 2. Main Results for Mixed Logic-Systems.}\par\medskip 

Let $\b L$ be any, at least, denumerable set that represents a language embedded into the standard superstructure $\cal M$ and ${\cal P}$ be the set-theoretic power set operator. Let $1 \leq n \in \nat,\ [1,n] =\{i\mid (i\in \nat)\land (1\leq i\leq n)\},$ and $ f\colon [1,n] \to \b L,\ g\colon [1,n] \to \Hyper {\b L},\ h\colon [1,n] \to \b L$ and for each $i,j\in [1,n], \ f(i)\not=h(j),\ g(i) \in \Hyper {\b L}- \b L.$ The functions $f,\ h$ can be  formed informally prior to embedding them into the standard superstructure $\cal M.$ Due to the finite domain and range of the functions $f,\ h$ and the identification process, it follows that we can consider $\hyper f = f, \Hyper h = h$ and, notationally, for each $i \in [1,n],\ f(i) = \b a_i, \ g(i) = \lambda_i,\ h(i) = \b b_i.$ Of course, the hyper-form can be retained. The ternary relation $R=\{(\b a_i,\lambda_i,\b b_i)\mid i \in [1,n]\}$ is called a {\it mixed logic-system}, and it is not a member of the carrier of the superstructure $\cal M.$ Let $\rm T_3 = \{ x\mid (x \in \power {L \times L \times L)} \land (\vert x \vert = n)\}.$ In Herrmann (2001, p. 94), it is shown that there is a map $\rm H_3\colon \r T_3 \to {\cal C}_{\r f},$ where ${\cal C}_{\r f}$ is the set of all finitary consequence operators defined on $\r L.$ For each $\r x \in \r T_3,$ the finite consequence operator $\r H_3(\r x)$ is defined by a logic-system using the rules of inference $\r x$ and the logic-system algorithm. For each $\rm X \subset L$, the set of all members of $\rm L$ deduced from $\rm X$ using these rules and the algorithm is $\r Y$ and this is denoted by $\rm X \vdash_x Y.$ Hence, $\rm X \vdash_x = H_3(x)(X)$  and each $\rm \vdash_x$ determines a binary relation that corresponds to that generated by $\rm H_3(x).$  \par 

In what follows, the notion of a finite subsequence includes the case where the domain is $\{1\}.$ In Herrmann (1993), there are three different fonts used for the mathematical expressions. Roman fonts refer to informal mathematics as expressed in general set-theory. Italic and bold face notation is used for statements about the objects in the standard and nonstandard superstructures. Further, $\omega$ represents the natural numbers in the informal general set-theory, while $\nat$ is isomorphic to $\omega$ with an additional set-theoretic restriction. In the published version of Herrmann (2001), and many others, the use of roman fonts for mathematical expressions is not allowed. For this reason, in these papers, $\nat$ replaces $\omega$ when informal mathematical expression are considered. \par\medskip

{\bf Theorem 2.1.} {\it Consider the defined ternary relation $R = \{(\b a_i,\lambda_i,\b b_i)\mid i \in [1,n]\},\ D_i= \{\b a_i, \lambda_i\},$ for each $i \in [1,n]$ and ${\cal C}_{\b f}$ be ${\cal C}_{\r f}$ embedded into the standard superstructure $\cal M.$ Then $R \in \Hyper {\b T_3}, \ \Hyper {\b H}_3(R)=C \in \Hyper {\cal C}_{\b f},$ $C$ is nonstandard and, for each $i\in [1,n],\ D_i$ is internal. Further, for each internal $X \in \Hyper (\power {\b L}),$ if $D_i\not\subset X,$ for each $i \in [1,n],$  then $C(X) =X$ and if $D_k \subset X$ for some $k \in [1,n],$ then there exists $m\in [1,n]$ and a finite subsequence $\{D_{i_j}\}$ of $\{D_i\}$ such that $D_k \subset X$ if and only if $k = i_j$ for some $j \in [1,m]$ and, for each $j \in [1,m]\},\ C(D_{i_j})= D_{i_j} \cup\{\b b_{i_j}\}$  and $C(X)= X \cup \{\b b_{i_j}\mid j\in [1,m]\}.$}  \pars
Proof. Consider any nonempty informal $\rm R' = \{(a_i, d_i, \r b_i) \mid i \in [1,n]\} \in \r T$ where, for each $\rm i,k\in [1,n],\ a_i \not= b_k,\ d_i\not= b_k.$ Let $\rm D_i= \{ a_i, d_i\}$ for each $\rm i \in [1,n].$  Let $\rm C\colon \power {\r L} \to \power {\r L}.$ We first establish that if $\rm H_3( R') = C$, then 
 for each $\rm X \in \power { L}$ (a) if $\rm D_i \not\subset X$ for each $\rm i \in [1,n],$  then $\rm C( X) =X$ and (b) if $\rm D_k\subset  X$ for some $\rm k \in [1,n],$ then there exists $\rm m\in [1,n]$ and a finite subsequence $\rm \{D_{i_j}\}$ of $\rm \{D_i\}$such that $\rm D_k \subset X$ if and only if $\rm k = i_j$ for some $\rm j \in [1,m]$ and, for each $\rm j \in [1,m],\ C( D_{i_j})= D_{i_j}\cup\{ b_{i_j}\}$  and $\rm C(X)= X \cup \{ b_{i_j}\mid j\in [1, m]\}.$ \par
We use the rules that generate a consequence operator from a logic-system (Herrmann, 2001, p. 94). First, there is a unique $\rm C \in {\cal C}_{ f}$ such that $\rm H_3( R') =  C.$ Let $\rm X \in \power { L}.$ (a) By the insertion rule $\rm X \subset  C(X).$ If $\rm D_i \not\subset X$ for each $\rm i \in [1,n]$, then by the coordinate rule it follows that $\rm X =  C(X).$ (b) If $\rm D_k \subset  X$ for some $\rm k \in [1,n],$ then the set $\rm M =\{k\mid  D_k \subset  X\}$ is nonempty, finite and orderable by the first element method. Let $\rm \vert M \vert = m.$ Then $\rm 1\leq m \leq n$ and there exists $\rm m\in [1,n]$ and a finite subsequence of sets $\rm \{D_{i_j}\}$ such that $\rm D_k \subset X$ if and only if $\rm k = i_j$ for some $\rm j \in [1,m].$  The coordinate rule yields $\rm C(D_{i_j}) =D_{i_j}\cup \{ b_{i_j}\}$ for each $\rm j\in [1,m].$ However, the value of $\rm C(X)$ is obtained only by insertion or (1) the coordinate rule applied to each $\rm D_i \subset X$ or (2) applied to members of $\rm \{ b_{i_j}\mid j \in [1,m]\}.$ The set $\rm \{ b_{i_j}\mid j\in [1,m]\}$ is all that can be deduced by (1) without considering (2). Since, for each $\rm i,k\in [1,n],\ a_i\not=b_k,\ d_i\not=b_k,$ the coordinate rule does not apply to any member of $\{\r b_{i_j}\mid j \in [1,m]\}.$ Hence, it follows that $\rm C(X) = X \cup \{b_{i_j}\mid j\in [1,m]\}.$\pars  
The defined partial sequences $f, g, h$ are internal from the Theorem 4.2.2 (ii) (Herrmann, 1991, p. 29) - The Internal Definition Principle. Also, from part (ii) of this theorem each member of the ternary relation $R$ is internal implies that $R,$ being a finite collection of internal objects, is internal.  This also follows directly from Theorem 4.2.2 part (ii).  Hence, $R \in \Hyper {\b T}_3$.  In like manner, each $D_{i}\in \Hyper (\power {\b L})$ for each $i \in [1,n].$ Moreover, since, in general, $\r H_3(\r R') \in {\cal C}_{\r f},$ for any $\r R' \in \r T_3,$ then $\Hyper {\b H}_3(R)=C \in \Hyper {\cal C}_{\b f},$ and $C$ is nonstandard, since $\b L$ is infinite. Note that, for each $i,j \in [1,n],\ \lambda_i \not= \b b_j.$ Further, $C$  satisfies, for each internal $X\in \Hyper (\power {\b L}),$ the appropriate formal *-transformed (a) and (b) statements. This task is left to the reader. This completes the proof. \qed\medskip

Now consider the mixed $g,\ f$ defined binary logic-system $R_1 = \{(\lambda_i,\b b_i)\mid i\in [1,n]\},$ $\r T_2 = \{\r x\mid (\r x \in \power {\r L \times \r L})\land (\vert \r x\vert = \r n)\},$ and $\r H_2 \colon \r T_2 \to {\cal C}_{\r f}.$ \par\medskip

{\bf Theorem 2.2.} {\it Consider the defined binary relation $R_1 = \{(\lambda_i,\b b_i)\mid i \in [1,n]\}.$ Let ${\cal C}_{\b f}$ be ${\cal C}_{\r f}$ embedded into the standard superstructure $\cal M.$ Then $R_1 \in \Hyper {\b B},\ \Hyper {\b H}_2(R_1)=C \in \Hyper {\cal C}_{\b f}$ and $C$ is nonstandard. Further, for each internal $X \in \Hyper (\power {\b L}),$ if $\lambda_i \notin X$ for each $i \in [1,n],$ then $C(X) =X$ and if $\lambda_k \in X$ for some $k \in [1,n],$ then there exists $m\in [1,n]$ and a finite subsequence $\{\lambda_{i_j}\}$ of $\{\lambda_i\}$ such that $\lambda_k \in X$ if and only if $k = i_j$ for some $j \in [1,m]$  and, for each $j \in [1,m]\},$ $C(\{\lambda_{i_j}\})= \{\lambda_{i_j}, \b b_{i_j}\}$  and $C(X)= X \cup \{\b b_{i_j}\mid j\in [1,m]\}.$}  \pars 
Proof. This is obtained by a simple modification of the proof of Theorem 1.\qed 
  
Of course, these two results can be extended to other mixed logic-systems.  Theorems 1 and 2 can be used to model certain hypotheses associated with the mind-brain problem. It is unusual to have rather simple set-theoretic statements characterize the correspondence between a logic-system and a consequence operator.\par\medskip

{\bf Theorem 2.3.} {\it Consider any nonempty $\rm R' = \{(a_i, d_i, \r b_i) \mid i \in [1,n]\} \in \r T_3,$ where for each $\rm i,k\in [1,n],\ a_i \not= b_k,\ d_i\not= b_k.$ Let $\rm D_i= \{ a_i, d_i\},$ for each $\rm i \in [1,n].$  Let $\rm C\colon \power {\r L} \to \power {\r L}.$ Then $\rm H( R') = C$ if and only if  
 for each $\rm X \in \power {L}$ ${\rm (}a{\rm )}$ if $\rm D_i \not\subset X,$ for each $\rm i \in [1,n],$  then $\rm C(X) =X,$ and ${\rm (}b{\rm )}$ if $\rm D_k\subset  X$ for some $\rm k \in [1,n],$ then there exists $\rm m\in [1,n]$ and a finite subsequence $\rm \{D_{i_j}\}$ of $\rm \{D_i\}$ such that $\rm D_k \subset X$ if and only if $\rm k = i_j$ for some $\rm j \in [1,m]$ and, for each $\rm j \in [1,m],\ C( D_{i_j})= D_{i_j}\cup\{b_{i_j}\}$ and $\rm C(X)= X \cup \{ b_{i_j}\mid j\in [1, m]\}.$} \par
Proof. $\Rightarrow$. The same proof as used for Theorem 1. \par
$\Leftarrow$ Suppose that you have the operator $\rm C'\colon \power {L} \to \power {L}$ and that $\rm C'$ satisfies (a) and (b). Let any $\rm X \in \power {L}.$ Suppose that $\rm D_i \not\subset X,$ for each $\rm i \in [1,n]$. Then $\rm C'(X) = X$ implies that $\rm X=C'(X)=C'(C'(X)).$  \par
Now suppose that $\rm D_k\subset  X$ for some $\rm k \in [1,n].$  Then there exists $\rm m\in [1,n]$ and a finite subsequence $\rm \{D_{i_j}\}$ of $\rm \{D_i\}$ such that $\rm D_k \subset X$ if and only if $\rm k = i_j$ for some $\rm j \in [1,m]$ and, for each $\rm j \in [1,m],\ C(D_{i_j})= D_{i_j}\cup \{ b_{i_j}\}$  and $\rm C(X)= X \cup \{ b_{i_j}\mid j\in [1,m]\}.$ Since for each $\rm i,k\in [1,n],\ a_i \not= b_k,\ d_i\not= b_k$ and $\rm D_i= \{ a_i, d_i\},$ for each $\rm i \in [1,n]$, then $\rm D_k \subset X \cup \{b_{i_j}\mid j\in [1,m]\}$ if and only if $\rm D_k \in \{\rm D_{i_j}\}.$ Hence, since the stated conditions apply to each member of $\power {\rm L}$, then $\rm C'(X \cup \{ b_{i_j}\mid j\in[1, m]\}) = (X \cup \{ b_{i_j}\mid j\in [1,m]\})\cup \{ b_{i_j}\mid j\in [1,m]\}=X \cup \{ b_{i_j}\mid j\in [1,m]\}= C'(X).$ Consequently, for the two cases, $\rm X \subset C'(X)) = C'(C'(X)).$\par
Let $\rm Z \in {\cal F}(X),$ where, in general, $\rm {\cal F}(Y)$ is the set of all finite subsets of $\rm Y.$ Since the conditions (a) and (b) apply to $\rm Z$, assume that $\rm D_i \not\subset Z,$ for each $\rm i \in [1,n].$ Then $\rm C'(Z) = Z.$ Now assume that $\rm D_{k} \subset Z$ for some $\rm k \in [1,n].$ Then there exists $\rm m_z\in [1,n]$ and a finite subsequence  $\rm \{D_{i_j}\}$ of $\rm \{D_i\}$ such that $\rm D_k \subset Z$ if and only if $\rm k = i_j$ for some $\rm j \in [1,m_z]$ and, for each $\rm j \in [1,m_z],\ C'( D_{i_j})=D_{i_j}\cup \{ b_{i_j}\}$  and $\rm C(X)= X \cup \{ b_{i_j}\mid j\in [1,m_z]\}.$ Clearly, since $\rm Z \subset X,\ m_z\leq m.$  Thus, $\rm C'(Z) = Z \cup \{b_{i_j}\mid j\in [1,m_z]\} \subset X \cup \{b_{i_j}\mid j\in [1,m]\}= C'(X).$ However, for each $\rm j \in [1,m], \  D_{i_j}$ is finite, $\rm D_{i_j} \subset X,$ and $\rm C'(D_{i_j}) = D_{i_j}\cup\{b_{i_j}\}.$ Considering, for each $\rm x \in X,\ Z = \{x\},$ then it follows that  $\rm C'(X) = \bigcup\{C'(Z)\mid Z \in {\cal F}(X)\}.$ These results imply that $\rm C' \in {\cal C}_f.$\par

Now consider $\rm H_3(R') = C$ and $\rm X \in \power{L}.$ Then from $\Rightarrow$ either the (a) or (b) hypotheses apply to $\rm X.$ But, the (a) and (b) hypotheses and the (b) statement ``there exists a $\rm m'\in [1,n]$ and a finite subsequence $\rm \{D_{i_j}\}$ of $\rm \{D_i\}$ such that $\rm D_k \subset X$ if and only if $\rm k = i_j$ for some $j \in [1,m']$'' only refer to the given sets $\rm D_i$ and $\rm X$ and not to the $\rm C.$ Thus, if (a) applies, then $\rm X = C(X) = C'(X)$. If (b) applies, then $\rm m' = m.$ Therefore, (a) and (b) imply that $\rm C(X) = C'(X)$ and this completes the proof. \qed

The following definitions identify a useful ordering for mixed logic-systems. Consider the same definitions for $f',\ g',\ h' $ and $\b a'_i,\ \lambda'_i, \ \b b'_i,$ respectively, but use $n'\in \nat,\ 1 \leq n'.$   The ternary relation $R=\{(\b a_k,\lambda_k,\b b_k)\mid k \in [1,n]\}$ is internal being a finite collection of internal objects.  Note that, due to the stated $f$ and $h$ properties, such logic-systems as $R$ and have been simplified for certain applications.  \par

By considering possible repetitions of the $\rm a_j$ and corresponding $\rm b_j$ as these members are embedded into ${\b L}$ and represented in the standard superstructure $\cal M$ via the partial sequences $f,\ h$ a special type of ordering can be defined. The following special ordering does not yield any mechanism but rather deals with a characteristic that can be considering as a type of ``weighting'' in the scientific sense. From a logic-system viewpoint, it yields a measure for an ``influencing process'' associated with perception. The idea is similar to the notion that ``repetition'' of a statement is a form of linguistic emphasis. \par\medskip
{\bf Definition 2.1.} Consider mixed logic-systems $R$ and $R' = \{(\b a'_k,\lambda'_k,\b b'_k)\mid k \in [1,n']\}$. Let $R_1 = \{(\b a_{j_k}, \lambda_{j_k},\b b_{j_k}) \mid k \in [1,m]\}\subset R$, where $j$ is the subsequence map defined on $[1,m]$ and $R_2 = \{(\b a'_{i_k}, \lambda'_{i_k},\b b'_{i_k}) \mid k \in [1,m']\}\subset R',$ and $i$ is the subsequence map defined on $[1,m'].$ For the projection maps $p_1,\ p_3,$ let $p_1(R_1) = \{\b a_{j_k}\},\ p_3(R_1) = \{\b a_{j_k}\},\ p_1(R_2) = \{\b a'_{j_k}\}$ and $p_3(R_2) = \{\b b'_{j_k}\}$ Then the {\it influencing process} is stronger for $\rm b'_{i_1}$ than for $\rm b_{j_1}$ if $m' > m.$ [Note. Obviously, the $\lambda_k,\ k \in [1,m]$ and the $\lambda'_k,\ k \in [1,m']$  are distinct.] \par\medskip

{\bf Definition 2.2.} Consider mixed logic-systems $R=\{(\lambda_k,\b b_k)\mid k \in [1,n]\},\  R'=\{(\lambda'_k,\b b'_k)\mid k \in [1,n']\}$. Let $R_1 = \{(\lambda_{j_k},\b b_{j_k}) \mid k \in [1,m]\}\subset R,$ where $j$ is the subsequence map defined on $[1,m].$ Let $R_2 = \{\lambda'_{i_k},\b b'_{i_k}) \mid k \in [1,m']\}\subset R',$ where $i$ is the subsequence map defined on $[1,m'].$   For the projection map $p_2,$ let $p_2(R_1) = \{\b b_{j_k}\}$ and $p_2(R_2) = \{\b b'_{j_k}\}.$ Then the {\it influencing process} is stronger for $\rm b'_{i_1}$ than for $\rm b_{j_1}$ if $m' > m.$ \par\medskip

In applications of definitions 2.1 and 2.2, the emphasis produced by the repeated members may or may not be an actual perceived ``repetition,'' from the viewpoint of the ultralogic operators. For example, let internal $X= \{\b a_{j_k}, \lambda_{j_k}\mid k \in [1,m]\}, \ X' = \{a'_{i_k},\lambda'_{i_k}\mid k \in [1,m']\}$ and the corresponding $R_1$ and $R_2$ satisfy the requirements of Definition 2.1. Then for ultralogics $C,\ C'$ of Theorem 2.1 and from its proof, it follows that $C(X)-X = \{\b b_{j_1}\},\ C'(X') -X = \{\b b'_{i_1}\}.$ Thus, if $m' > m,$ then the stronger influencing process may be more relative to ``how'' $\rm b'_{i_1}$ is perceived. From the viewpoint of physical processes, this can be interpreted as stating that, as perceived, the process that yields $\rm b'_{i_1}$ is stronger than the process that yields $\rm b_{j_1}$ even if the process itself is not known. \par\medskip

\noindent {\bf 3. Hyperfinite Logic-Systems.}\par\medskip 

Let ${\cal F}'$ denote the finite power set operator, where the range of ${\cal F}'$ does not contain the empty set. Let $\emptyset \not=X_i \subset \b L_i = \b L,$ where $i \in [1,n] \subset \nat, \ 1 < n.$ Let $B = X_1 \times \cdots \times X_n \subset \b L_1 \times \cdots \times \b L_n$. Then each $X_i$ is finite if and only if $B$ is finite. For nonempty internal $Y_i \subset \Hyper {\b L}, \ Y_1 \times \cdots \times Y_n, \ (i \in [1,n] \subset \nat, \ 1 < n)$ is an internal subset of $\Hyper {\b L}_1 \times \cdots \times \Hyper {\b L}_n$ by the Theorem 4.2.2 (ii) in Herrmann (1991, p. 29). Let nonempty hyperfinite $Z_i \subset \Hyper {\b L}, \ i \in [1,n], \ 1 < n \in  \nat \subset \hypernat.$ Since each such $Z_i \in \hyper {\cal F}'(\Hyper {\b L}),$ then each such $Z_i$ is internal. Let $\rm T_n' = \{ x \in \mid x \in {\cal F}'(L_1 \times \cdots \times L_n)\}.$ Then once again, from Herrmann (2001, p. 94), there is a function $\rm H_n'\colon T_n' \to {\cal C}_f$. By *-transfer, we have\par 
\medskip

{\bf Theorem 3.1.} {\it For each hyperfinite $A \in \hyper {\cal F}'(\Hyper {\b L_1} \times \cdots \times \Hyper {\b L_n}),$ there exists $C_A \in \Hyper {\cal C}_f$ such that $\Hyper {\b H_n'}(A) = C_A.$ Further, for each $X \in \Hyper {\power {\Hyper {\b L}}}, \ X {\vdash_A} = C_A(X).$} \qed 

{\bf Example 3.1} Let nonempty $d \in \Hyper {\power{\Hyper {\b L}}}.$ Then $d \in \Hyper {\cal F}'(\Hyper {\b L})$ if and only if there is a $\lambda \in \hypernat$ and an internal bijection $f\colon [0,\lambda] \to d.$ Let $\lambda > 1$ and $R = \{(x,y)\mid \forall i((i \in [0,\lambda -1])\land  (x = f(i))\land (y = f(i+1))\}.$ For any $g \in d$, $d -\{g\}$ is hyperfinite by *-transfer from the finite case. Thus, the range and domain of $R$ are hyperfinite. Consequently, $R \in \Hyper {\cal F}'(\Hyper {\b L} \times \Hyper {\b L}).$ Let $\Hyper {\b H}'_2(R)= C \in \Hyper {\cal C}_f.$ Then $C(\{f(0)\})= d$ and, for $i \not= 0,\  C(\{f(i)\}) \not= d.$ \bigskip

\centerline{\bf References}\par\medskip
\id{G}eiser, James R. (1968). ``Nonstandard logics,'' {\it J. of Symbolic Logic}, 33(2):236-250. \pars
\id{H}errmann, Robert A. (2006). ``General logic-systems and finite consequence operators,'' {\it Logica Universalis}, 1:201-208.\hfil\break http://arxiv.org/abs/math/0512559\pars  
\id{H}errmann, Robert A. (2001). ``Hyperfinite and Standard Unifications for Physical Theories,'' {\it Intern. J. Math. Math. Sci.}, 28(2):93-102.\hfil\break http://arxiv.org/abs/physics/0105012\pars

\id{H}errmann, Robert A. (1993). {\it The Theory of Ultralogics}, \hfil\break http://arxiv.org/abs/math.GM/9903081 
http://arxiv.org/abs/math.GM/9903082\pars

\id{H}errmann, Robert A. (1991). {\it Nonstandard Analysis Applied to Advanced Undergraduate Mathematics - Infinitesimal Modeling},\hfil\break http://arxiv.org/abs/math.GM/0312432\pars 

\id{R}obinson, Abraham (1963). ``On languages which are based on non-standard arithmetic,'' {\it Nagoya Math. J.}, 22:83-117.\pars

\end